\documentclass[a4paper, 
11pt
]{article}

\usepackage{subfig}
\usepackage{pgfplots}
\usepackage{pgfplotstable}
\usepackage{mathtools}
\usepackage{multicol}
\usepackage{comment}
\usepackage{booktabs}
\pgfplotsset{compat=1.5}
\usepackage{amssymb}
\usepackage{url}

\usepackage{pdfsync}
\usepackage{float}
\usepackage{tabularx}
\usepackage{enumerate}
\usepackage{array}
\usepackage{xspace}

\usepackage{authblk}
\usepackage{symbols}

\newcommand{\mupar}{\ensuremath{\boldsymbol{\mu}}}
\newcommand{\etapar}{\ensuremath{\boldsymbol{\eta}}}

\begin{document}

\title{An efficient shape parametrisation by free-form deformation enhanced
by active subspace for hull hydrodynamic
ship design problems in open source environment}

\author[]{Nicola~Demo\footnote{nicola.demo@sissa.it}}
\author[]{Marco~Tezzele\footnote{marco.tezzele@sissa.it}}
\author[]{Andrea~Mola\footnote{andrea.mola@sissa.it}}
\author[]{Gianluigi~Rozza\footnote{gianluigi.rozza@sissa.it}}

\affil[]{Mathematics Area, mathLab, SISSA, via Bonomea 265, I-34136 Trieste, Italy}

\maketitle

\begin{abstract}
In this contribution, we present the results of the application of a parameter
space reduction methodology based on active subspaces property to the hull hydrodynamic
design problem. In the framework of such typical naval architecture problem, several
parametric deformations of an initial hull shape are considered to assess the influence
of the shape parameters considered on the hull total drag. The hull resistance, which
is the performance parameter associated with each parametric hull, is typically
computed by means of numerical simulations of the hydrodynamic flow past the ship.
Such problem is extremely relevant at the preliminary stages of the ship design, when
several flow simulations are typically carried out by the engineers to establish a
certain sensibility on the total drag dependence on the hull geometrical parameters
considered and on other physical parameters.
Given the high number of geometric and physical parameters involved –which might
result in a high number of time consuming hydrodynamic simulations– assessing whether
the parameters space can be reduced would lead to considerable computational cost
reduction at the design stage.
Thus, the main idea of this work is to employ the active subspaces to identify possible
lower dimensional structures in the parameter space, or to verify the
parameter distribution in the position of the control points. To this end, a fully automated
procedure has been implemented to produce several small shape perturbations of an
original hull CAD geometry which are then used to carry out high-fidelity flow
simulations in different cruise conditions and collect data for the active subspaces
analysis. To achieve full automation of the open source pipeline described, both the
free form deformation methodology employed for the hull perturbations and the high
fidelity solver based on unsteady potential flow theory, with fully nonlinear free
surface treatment, are directly interfaced with CAD data structures and operate using
IGES vendor-neutral file formats as input files. The computational cost of the fluid
dynamic simulations is further reduced through the application of dynamic mode
decomposition to reconstruct the final, steady state total drag value given only
few initial snapshots of the simulation. The active subspaces analysis is here applied
to the geometry of the DTMB-5415 naval combatant hull, which is a common benchmark in
ship hydrodynamics simulations, within the SISSA mathLab applied mathematics lab. The
contribution will discuss several details of the implementation of the tools developed,
as well as the results of their application to the target engineering problem.
\end{abstract}

\section{Introduction}
\label{sec:intro}
Active subspaces (AS) has emerged as a useful parameter study technique in
the last years (see~\cite{constantine2015active}). It allows to identify a low
dimensional structure in the parameter space of the scalar output
function of interest. It uses a linear combination of all the
parameters as new variable in order to reveal such lower dimensional
behaviour of the output function. The aim of the present work is
to investigate on the application of AS to reduce the number of parameters
involved in the deformation of a particular hull geometry. Not only
such investigation might result in a direct indication of a model reduction
strategy for the problem at hand, but it can also illustrate how the
AS algorithm can be applied to generic design problems to check
whether a further reduction of dimensionality is possible.

To feed the AS algorithm, the relationship between the shape (input) parameters
and the different fluid dynamic performance (output) parameters
must be assessed. This has been done through a campaign of high-fidelity fluid
dynamic simulations. In each of such simulations, the hull model has been
morphed according to the deformation corresponding to the sample point in the parameters
space.

As for the shape parametrization strategy, in the framework of reduced
order modeling (ROM) we adopted a so-called general purpose
shape morphing, that is a technique that does not depend on the
particular shape to deform. As possible choices we mention free form
deformation (FFD), inverse distance weighting (IDW) or radial basis
functions (RBF) interpolation. All these methods basically involve the displacement of
some control points in order to deform a domain. The displacements of
these control points are what we identify as parameters. In this work
we use the FFD since it allows global deformations of the geometry
with only few parameters (see~\cite{sederbergparry1986}).

The fluid dynamic solver used in this work is based on fully nonlinear potential
flow model, implemented in the software WaveBEM (\cite{molaEtAl2013,MolaHeltaiDeSimone2017}). Being
directly interfaced with CAD data structures, the solver is able to import the
deformed geometries, automatically generate the computational grids and carry out
the simulation with no need for human interaction. The unsteady flow simulations
were significantly accelerated through the use of dynamic mode decomposition (DMD)
methodology. Since its introduction in the fluid mechanics community, DMD has
emerged as a useful tool for analyzing the
dynamics of nonlinear systems, by the approximation of the Koopman
operator (see~\cite{schmid2010dynamic}). We will exploit the linear finite dimensional Koopman
operator to reduce the computational cost of a single high-fidelity
simulation. The DMD method allows to reconstruct the system dynamics
through few snapshots describing the state of the system. With this
approximation we can also infer on the future evolution of the system,
allowing us to simulate only few seconds using the high-fidelity
solver. We underline that both the pressure and velocity fields are
reconstructed through the DMD.

The paper will briefly present all the aforementioned methodologies, and
details of their application to the hull design problem considered.
The results of the combined application of DMD and AS to the outputs of the fluid dynamic
simulations will also be presented and discussed.

\section{Estimation of the resistance  of a hull advancing in calm water}
\begin{figure}[htb]
\centering
\includegraphics[width=0.9\textwidth]{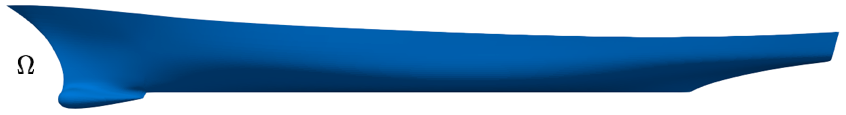}
\caption{Reference domain $\Omega$, that is the DTMB 5415 hull.}
\label{fig:original_hull}
\end{figure}

In this section we introduce the problem of the estimation of the
resistance of a ship advancing in calm water. The hull shape considered in
this work is that of the DTMB 5415, which was conceived for the preliminary design of
a US Navy Combatant ship and includes a sonar dome and a transom stern (see
Figure~\ref{fig:original_hull}). Given the abundance of experimental data
available in the literature (among others, we cite~\cite{olivieri2001towing})
has become a common benchmark for naval hydrodynamics simulation tools.

We denote with $\Omega \subset \mathbb{R}^3$ the domain (see Figure
\ref{fig:original_hull}) associated with our model hull.  
More specifically, $\Omega$ is our \emph{reference} domain, and corresponds
with the undeformed DTMB 5415 hull --- the latter assumption is not
fundamental for the remainder of the paper, but is convenient for practical reasons.
We must here remark that the domain considered in the fluid dynamic simulations
is in principle the volume of water $\Omega_w$ surrounding the hull. Further details
about the fluid dynamic domain will be provided in the following sections.

We define the shape morphing $\mathcal{M}(\boldsymbol{x}; \boldsymbol{\mu}): \mathbb{R}^3 \to \mathbb{R}^3$
that maps the reference domain $\Omega$ into the deformed
domain $\Omega(\boldsymbol{\mu})$, namely:
\begin{equation*}
\Omega(\boldsymbol{\mu}) = \mathcal{M}(\Omega; \boldsymbol{\mu}).
\end{equation*}
It is quite natural to infer that the flow field, and thus the result of the fluid dynamic simulations, will
depend on the specific hull shape considered. In turn, such shape is associated to the parameters defining the
morphing $\mathcal{M}$ --- which will be extensively  defined in the next sections. One of the main purposes
of this contribution is then to investigate the effect of the morphing parameters on the total resistance,
the main fluid dynamic performance parameter.
We must here remark that one of the geometrical quantities having the most effect on the resistance is
the immersed volume of a hull shape, as higher volumes will generate higher drag values.
If designers only consider the absolute resistance value, they might then disregard hull shapes
that have relatively good drag performances despite having higher immersed volume. This is of course
undesired and could be avoided adopting different strategies. A first possibility is considering a
resistance value which is nondimensionalized using a measure of the displaced volume. In alternative,
the morphing parameters could be constrained so as to impose a fixed hull immersed volume.
Unfortunately, the former solution poses some problem in the identification
of the most suitable nondimensionalization strategy resulting in resistance output indices truly independent of the
hull immersed volume. As for the latter possibility, it would lead to an undesired complication of the shape
parametrization methodology. In this work the effect of immersed volume variations on the hull resistance is
taken into account in a more natural fashion. In the fluid dynamic simulations, the rigid motions of the hull are also
considered. In this way, once the ship displacement is imposed, the hull shape considered in the simulation
reaches the equilibrium position under the action of gravity and hydrodynamic forces. When a shape morphing
is characterized by a higher volume, the computed sink will reduce to obtain the same vertical
component of the hydrodynamic force. This ensures that each shape is compared on a level ground.

By a practical standpoint, once a point in the parameter domain $\mathbb{D}$ is identified, the specific hull
geometry is provided to the fluid dynamic solver, which
carries out a flow simulation to come up with a resistance estimate. In this framework free form deformation has been
employed for the generation of a very large number of hull geometries based on the DTMB 5415
naval combatant hull shape morphing. Each geometry generated has been used to set up a high-fidelity
hydrodynamic simulation with the desired ship displacement and hull speed. Since the serial time
dependent high-fidelity simulations are quite time consuming, taking approximatively 24 hours to
reach a steady state solution, we adopted a reduction strategy based on dynamic mode decomposition (DMD)
to cut the overall computational cost of each simulation to roughly 10 hours.
The output resistances for all the configurations tested have been
finally analyzed by means of active subspaces (AS) in order to verify
if a further reduction in the parameter space is feasible, which could
significantly speed up the work of designers.

In the next sections, we will provide a brief description of the
unsteady fully nonlinear potential fluid dynamic model used to carry out the
high-fidelity simulations. We refer the interested
reader to~\cite{molaEtAl2013,mola2016ship,MolaHeltaiDeSimone2017} for
further information on the fully nonlinear potential free surface model,
on its application to complex hull geometries, and on the treatment of
the hull rigid motions respectively. In addition, we will describe the free form
deformation morphing methodology implemented, and describe the dynamic
mode decomposition strategy used to cut the
computational time of the simulations. Finally, the active subspaces used to
analyze the resistance dependence on the morphing parameters will be presented.

\section{Fully nonlinear potential model}

In the simulations we are only considering the motion of a ship advancing
at constant speed in calm water. For such reason we solve the problem
in a \emph{global}, translating reference frame $\widehat{XYZ}$,
which is moving with the constant horizontal velocity of the
boat $\Vb_\infty = (V_\infty,0,0)$. Thus, the $X$ axis
of the reference frame is aligned with $\Vb_\infty$, the $Z$ axis is
directed vertically (positive upwards), while the $Y$ axis is directed
laterally (positive port side). 

As aforementioned, the domain $\Omega_w(t)$ in which we are interested in
computing the fluid velocity and pressure is represented by the portion of
water surrounding the ship hull. The time varying shape of such
domain --- and in particular that of its boundary with the air above --- is
one of the unknowns of the fluid dynamic problem. By convention, we place the
origin of the vertical axis $Z$ in correspondence with the undisturbed free surface
level, and we start each simulation at time $t=0$ from such undisturbed
configuration. Thus, $\Omega_w(t=0) = \mathbb{R}^3_{Z-}\ \Omega$, where
$\mathbb{R}^3_{Z-}$ indicates the lower half-space of $\mathbb{R}^3$, for which $Z\leq 0$.

If overturning waves are not observed (which is typically the case for low cruise velocities
typical of a ship), the domain $\Omega_w(t)$ is simply connected. So under the assumptions of
irrotational flow and non viscous fluid the velocity field $\vb(\Xb,t)$
admits a representation through a scalar potential function
$\Phi(\Xb,t)$, namely
\begin{equation}
\label{eq:potential-definition}
\vb = \nablab\Phi = \nablab\left(\Vb_\infty\cdot\Xb + \phi \right)\qquad\qquad
\forall \ \Xb\in \Omega_w(t),
\end{equation}
in which $\phi(\Xb,t)$ is the \emph{perturbation potential}.
Under the present assumptions, the equations of motion simplify to the unsteady
Bernoulli equation and to the Laplace equation for the perturbation
potential:
\begin{subequations}
  \label{eq:incompressible-euler-potential}
  \begin{alignat}{3}
    \label{eq:bernoulli}
    & \frac{\Pa\phi}{\Pa t}+\frac{1}{2}\left|\nablab\phi+\Vb_\infty\right|^2
      +\frac{p-p_a}{\rho}-\gb\cdot\Xb
       = C(t) \qquad & \text{ in } \Omega_w(t) ,\\
    \label{eq:incompressibility-potential}
    & \Delta \phi = 0 & \text{ in } \Omega_w(t) ,
  \end{alignat}
\end{subequations} 
where $C(t)$ is an arbitrary function of time, and $\gb = (0,0,-g)$,
is the gravity acceleration vector, directed along the $z$ axis. In the
equation the constant water density $\rho$ and atmospheric reference pressure
$p_a$ are also appearing.
A noteworthy characteristic of the fluid dynamic problem at hand, is that the
unknowns of such mathematical problem $\phi$ and $p$ are
uncoupled. This means that the solution of the Poisson problem in Eq.~\eqref{eq:incompressibility-potential}
can be obtained independently of the pressure field. Once such solution
is obtained, the pressure can be obtained through a  postprocessing step
based on Bernoulli Eq.~\eqref{eq:bernoulli}. Thus, the Laplace equation is the
governing equation of our model. Such equation is complemented by non penetration
boundary conditions on the hull surface $\Gamma^{b}(t)$ and water basin bottom boundary
$\Gamma^{bot}(t)$, and by homogeneous Neumann boundary conditions on the truncation
boundaries $\Gamma^{far}(t)$ of the numerical domain. The bottom of the basin is located
at a depth corresponding to 2 boat lenghts, while the truncation boundaries are
located approximatively at a distance from the boat of 6 boat lengths in the longitudinal
direction $X$ and of 2 boat lengths in the lateral direction $Y$). 
On the water free surface $\Gamma^{w}(t)$, we employ the kinematic and dynamic Semi-Lagrangian fully nonlinear boundary
conditions, which respectively read
\begin{eqnarray}
\frac{\delta\eta}{\delta t} &=& \frac{\Pa\phi}{\Pa z}+\nablab\eta\cdot\left(\wb -\nablab\phi-\Vb_\infty\right)
\qquad \text{ in } \Gamma^{w}(t) ,
\label{eq:fsKinematicBeck}\\
\frac{\delta\phi}{\delta t} &=& -g\eta +
\frac{1}{2}|\nablab\phi|^2 + \nablab\phi\cdot\left(\wb -\nablab\phi-\Vb_\infty\right)
\qquad \text{ in } \Gamma^{w}(t).
\label{eq:fsDynamicBeck}
\end{eqnarray}  
The former equation expresses the fact that a material point moving on the
free surface will stay on the free surface --- here assumed to be
a single valued function $\eta(X,Y,t)$ of the horizontal
coordinates $X$ and $Y$. The latter condition
is derived from Bernoulli Eq.~\eqref{eq:bernoulli},
under the assumption of constant atmospheric pressure on the water surface.
This peculiar form of the fully nonlinear boundary
conditions was proposed in~\cite{beck1994}. Eq.~\eqref{eq:fsKinematicBeck}
allows for the computation of the vertical velocity of markers which move on the
water free surface with a prescribed horizontal speed $(w_X,w_Y)$.
Eq.~\eqref{eq:fsDynamicBeck} is used to obtain the velocity potential
values in correspondence with such markers. The resulting
vector $\wb=(w_X,w_Y,\frac{\delta\eta}{\delta t}) = \dot{\Xb}$ is the time
derivative of the position of the free surface markers.
In this work, such free surface markers are chosen as the free surface nodes
of the computational grid. To avoid an undesirable mesh nodes drift
along the water stream, the markers arbitrary horizontal velocity is
set to 0 along the $X$ direction. The $Y$ component of the water nodes in
contact with the ship --- which is moved according with the computed linear and
angular displacements --- is chosen so as to keep such nodes on the hull surface.
As for the remaining water nodes, the lateral velocity value is set  
to preserve mesh quality.

\subsection{Three dimensional hull rigid motions}

The equations governing the fluid dynamic problems are interfaced with
the equations describing the dynamics of the hull --- here supposed rigid.
The evolution equation for the hull baricenter position $\Xb^G(t)$ is
obtained via the linear momentum conservation equation, which in the case
of our hydrodynamics simulation framework reads
\begin{equation}
\label{eqHullLinCons}
m_s\ddot{\Xb}^G(t) = m_s\gb + \Fb^w(t).
\end{equation}
In Eq.~\eqref{eqHullLinCons}, $m_s$ represents the mass of the ship, while the
hydrodynamic force vector $\Fb^w(t)$ is obtained as
the sum of the pressure and viscous forces on the hull. 

The angular momentum conservation is instead written to obtain an evolution
equation for the angular velocity $\omegab$ of the hull, namely
\begin{equation}
\label{eqHullAngCons}
R(t) I^G R(t)^T\dot{\omegab}(t) + \omegab(t) \times R(t) I^G R(t)^T \omegab(t) = \Mb^w(t),
\end{equation}

where the hydrodynamic moment vector $\Mb^w(t)$ is the sum of the
moment about the ship center of gravity of the pressure and viscous forces
on hull, propeller and appendages. At a closer look, a further unknown
appears in Eq.~\eqref{eqHullAngCons}. The rotation matrix $R(t)$ is in fact
to be determined during the simulation. Such matrix allows for the conversion
of points coordinates from the global reference frame $\widehat{XYZ}$ to
the hull attached reference frame $\widehat{xyz}$ in which the hull geometry
is described and the hull moment of inertia $I^G$ is computed, namely  
$\Xb(t) = R(t)\xb+\Xb^G(t)$.
To obtain a set of equations which govern the time evolution of $R(t)$, we
resort to hull quaternions (again, we refer to~\cite{shoemake1985} and~\cite{mola2016ship} for further
detail on quaternions and their application to hull dynamics respectively). For the present
work purposes a quaternion is a  particular type of four element vector, namely
$\qb(t) =  \left[s,\vb\right] = s + v_X(t)\eb_X + v_Y(t)\eb_Y + v_Z(t)\eb_Z $.
If the norm of a quaternion (defined as $||\qb(t)|| = \sqrt{s(t)^2+v_X(t)^2+v_Y(t)^2+v_Z(t)^2}$)
has a unit value, the quaternion identifies a rotation matrix by means
of the relation

\begin{equation}
\label{eqQuaterToRotMat}
R =
\left[\begin{array}{c c c}
1-2v_Y^2-2v_Z^2 &  2v_Xv_Y-2s v_Z & 2v_Xv_Z+2s v_Y \\
2v_Xv_Y+2s v_Z & 1-2v_Y^2-2v_Z^2 & 2v_Yv_Z-2s v_X \\
2v_Xv_Z-2s v_Y & 2v_Yv_Z+2s v_X & 1-2v_Y^2-2v_Z^2
\end{array}\right],
\end{equation} 
in which to lighten the notation we omitted the time dependence of both $R(t)$ and
the components of $\qb(t)$. Thus, a unit quaternion is associated with the hull, and
allows for the computation of the rotation matrix at all times. The equation describing
the time evolution for the hull unit quaternion $\qb(t)$ is 
\begin{equation}
\label{eqQuaterEvo}
\dot{\qb}(t) = \dsfrac{1}{2}\omegab_q(t)\qb(t),
\end{equation}
where $\omegab_q(t) = \left[0,\omegab(t)\right]$ is the quaternion associated with
the angular velocity vector $\omegab(t)$, and the quaternion internal product is defined as
\begin{equation}
\label{eqQuatProd}
\qb_1\qb_2 = \left[s_1,\vb_1\right]\left[s_2,\vb_2\right] =
\left[s_1s_2-\vb_1\cdot\vb_2\, , \, s_1\vb2+s_2\vb1+\vb_1\times\vb_2\right].
\end{equation}
As quaternions only have four entries,
there only is one extra variable used to describe the three degrees of freedom of a
three dimensional rotation. This is far less redundant than writing an evolution equation for each
of the nine entries in a three dimensional rotation matrix. The residual redundant degree
of freedom involved in the hull quaternion evolution is easily eliminated imposing a unit
constraint on its norm, so as to ensure that the hull quaternion
always represents a pure rotation.

\section{Discretization and numerical solution }

The boundary value problem described in the previous section is governed by the linear Laplace operator.
Nonetheless, the presence of boundary conditions
Eqs.~\eqref{eq:fsKinematicBeck} and \eqref{eq:fsDynamicBeck} make the problem nonlinear. More
sources of nonlinearity are given by continuous change of the domain shape over time
and by the arbitrary shape of the ship hull. At each time instant, we will
compute the values of the unknown potential and node displacement fields by
solving a specific nonlinear problem resulting from the spatial and time discretization
of the original boundary value problem. The spatial discretization of the Laplace
problem is based upon a boundary integral formulation described in
detail in~\cite{molaEtAl2013} and~\cite{giulianiEtAl2015}. The domain
boundary is partitioned into quadrilateral cells, on
which bi-linear shape functions are used to approximate the surface, the flow
potential values, and the normal component of its surface gradient. The resulting iso-parametric
Boundary Element Method (BEM, see~\cite{brebbia}) consists in collocating a Boundary
Integral Equation (BIE) in correspondence with each node of the numerical grid, and
computing a numerical approximation the integrals appearing in such equation.
The linear algebraic equations obtained from such spatial discretization step are combined
with the Ordinary Differential Equations (ODE) derived from the finite element spatial 
discretization of the fully nonlinear free surface boundary conditions in Eqs.~\eqref{eq:fsKinematicBeck} and
\eqref{eq:fsDynamicBeck}. The spatial discretization described is carried out making use of the
deal.II open source library for C++ implementation of finite element
discretizations (see~\cite{BangerthHartmannKanschat2007,BangerthHeisterHeltai-2016-b}). 
To enforce strong coupling between the fluid and structural problem,
the aforementioned system of Differential Algebraic Equations (DAE) is complemented by the equations
of the rigid hull dynamics.
The fully coupled DAE system solution is time integrated by means of an arbitrary
order and arbitrary time step  implicit Backward Difference Formula (BDF) scheme implemented in the
IDA package of the open source C++ library SUNDIALS described in~\cite{hindmarsh2005sundials}. The potential
flow model described has been implemented in a stand alone C++ software, the main features of which are
described in~\cite{molaEtAl2013}.
 
As for the output of the simulations considered in this work, at each time step of the simulation,
the wave resistance is computed as
\begin{equation}
\label{eqLinForce}
R^w = \int_{\Gamma^b}p\nb\,d\Gamma\cdot\eb_X,
\end{equation}
where $p$ is the pressure value obtained introducing the computed potential in
Eq.~\eqref{eq:bernoulli}. The inviscid fluid dynamic model drag prediction
is then corrected by adding a viscous drag contribution obtained by means of the
ITTC-57 formula reported in~\cite{morrall19701957}. A full assessment of the accuracy of
the high-fidelity fluid structure interaction solver described is clearly beyond the scope
of the present work. Yet, for all the Froude numbers tested the computed total drag difference
with respect to the measurements reported by~\cite{olivieri2001towing} is less then 6\% (again,
we refer the interested reader to~\cite{molaEtAl2013,mola2016ship,MolaHeltaiDeSimone2017} for more details).
Given the fact that all the geometries tested are deformations of the present hull, it
is reasonable to infer that for each simulation carried out the accuracy of the high-fidelity model
prediction will be similar to that of the results presented.

No parallel version of the fluid dynamic solver has currently been developed. For such reason,
the computational time of 24 hours for each simulation (accounting for approximatively 5000 cells) is
not particularly competitive with respect to parallel RANS solvers. Yet, the extremely small size
of the simulation outputs and the simplicity of the surface mesh required for the simulations make the
present solver ideal for fully automated simulation campaigns based on shape parametrization.
Moreover, the solver is complemented with a mesh module directly interfaced
with CAD data structures based on the methodologies for surface mesh generation
described in~\cite{dassi2014curvature}. At the start of each simulation, such feature
allows for fully automated mesh generation from each hull shape assigned in the form of
a --- possibly non water-tight --- IGES geometry. This makes the fully nonlinear
potential solver described a particularly useful tool when used in the framework of
design pipelines which involve testing a high number of parametrised shapes. As we will
see, not only human user interaction is not required to carry out several hundreds simulations;
but in addition the extremely reduced dimension of the BEM solution files allows
for very simple file storing and for a fast post processing phase.

\section{Shape parametrization through free form deformation}
\label{sec:ffd}
In this section we will briefly describe how the free form deformation (FFD)
shape morphing strategy has been used to produce a parametrised set
of modified DTMB 5415 hulls to be  used in the fluid dynamic simulations.
FFD is a versatile parametrization technique used for shape optimization in
a variety of fields such as aerospace engineering, structural mechanics, and
biomedical engineering among others. For further insight on
the original formulation of FFD we suggest to
see~\cite{sederbergparry1986}, while for a recent work we suggest~\cite{salmoiraghi2018free}.
One of the main features of FFD, is that it does not directly manipulate the
geometrical object at hand. Instead it deforms a lattice of points
built around the object itself, manipulating the whole space in which
the geometry is embedded. This lattice has the topology of a hypercube
of dimension equal to the dimension of the
geometry we want to morph (3D in this work). The lattice is deformed
using a trivariate tensor-product of B-spline functions. This produces a
continuous and smooth deformation of the geometry. Such methodology has been
implemented into a stand alone Python package PyGeM~\cite{pygem}, which has been used
to obtain the morphed geometries considered in the present work.

\begin{figure}[htb]
\centering
    \includegraphics[trim={0 2cm 0 0},clip,width=.6\textwidth]{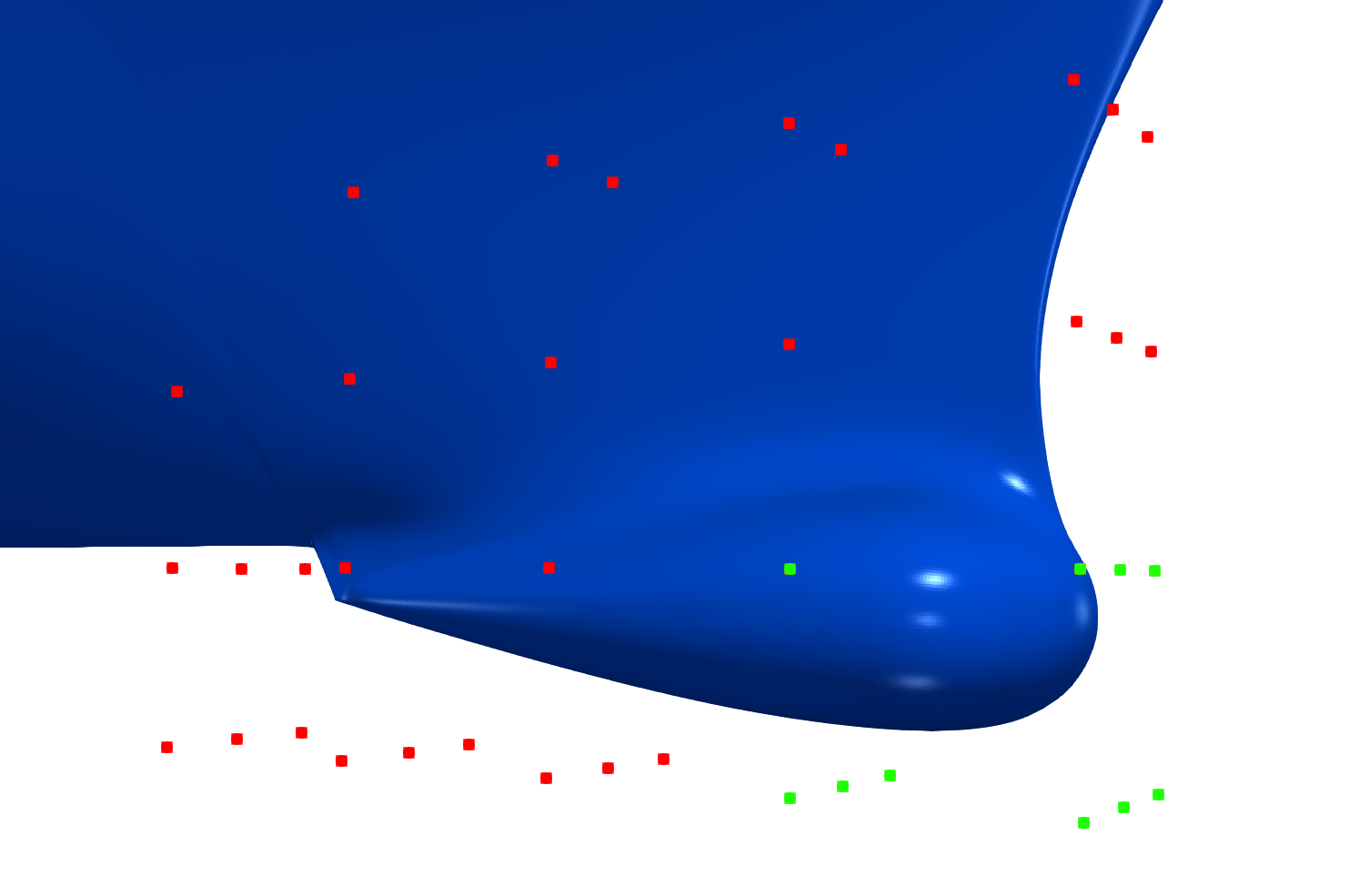}
\caption{Control points of the FFD lattice around the bulbous bow. The points
         which are displaced to generate the hull deformations are indicated with the green color.}
\label{fig:ffd_points}
\end{figure}

Figure~\ref{fig:ffd_points} displays a detail of the DTMB 5415 sonar dome,
surrounded by the points composing the FFD lattice. FFD is devised to only
deform objects which fall within the FFD lattice, so the picture indicates that
in this work we are only modifying the shape of the DTMB 5415 sonar dome.
The green lattice points in the picture are the ones that have
been moved to produce the deformations of domain $\Omega$ used in this work. 

The FFD procedure can be in fact subdivided into three steps. First, the physical
domain $\Omega$ is mapped to the reference domain $\widehat{\Omega}$
through a map~$\psi$. Then, some control points (the green ones in the picture, in our case)
$\boldsymbol{P}$ of the lattice are moved. Note that the displacement of such points are the parameters
$\boldsymbol{\mu}$ that identify a particular deformed geometry. Once the lattice has
been deformed, trivariate tensor-products of B-spline functions are used to compute the
map $\widehat{T}$ that associates the original lattice to the deformed one.
Finally, the back mapping from the deformed reference domain is applied to the deformed
physical domain $\Omega(\boldsymbol{\mu})$ using the map $\psi^{-1}$.
So it possible to express the FFD map $\mathcal{M}(\boldsymbol{x};
\boldsymbol{\mu}): \Omega \subset \mathbb{R}^3 \to
\Omega(\boldsymbol{\mu}) \subset \mathbb{R}^3$ by the composition
of these three maps, i.e.
\begin{equation}
\mathcal{M} (\cdot \, ; \boldsymbol{\mu}) = (\psi^{-1} \circ \widehat{T} \circ \psi)
(\cdot \, ; \boldsymbol{\mu}).
\end{equation}

\begin{figure}[htb]
\centering
\includegraphics[trim={0 3cm 0 0},clip,width=.44\textwidth]{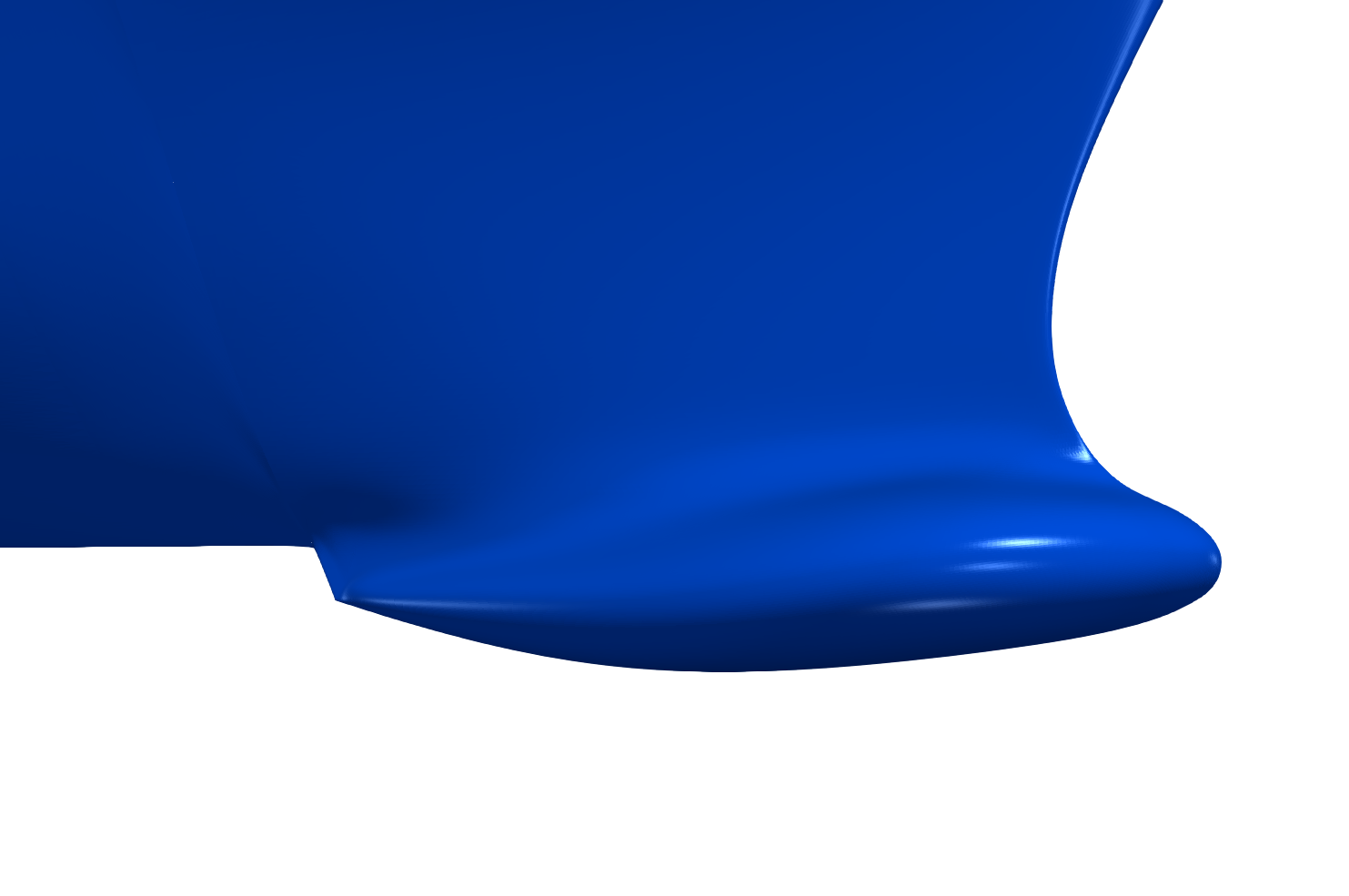}
\includegraphics[trim={0 3cm 0 0},clip,width=.44\textwidth]{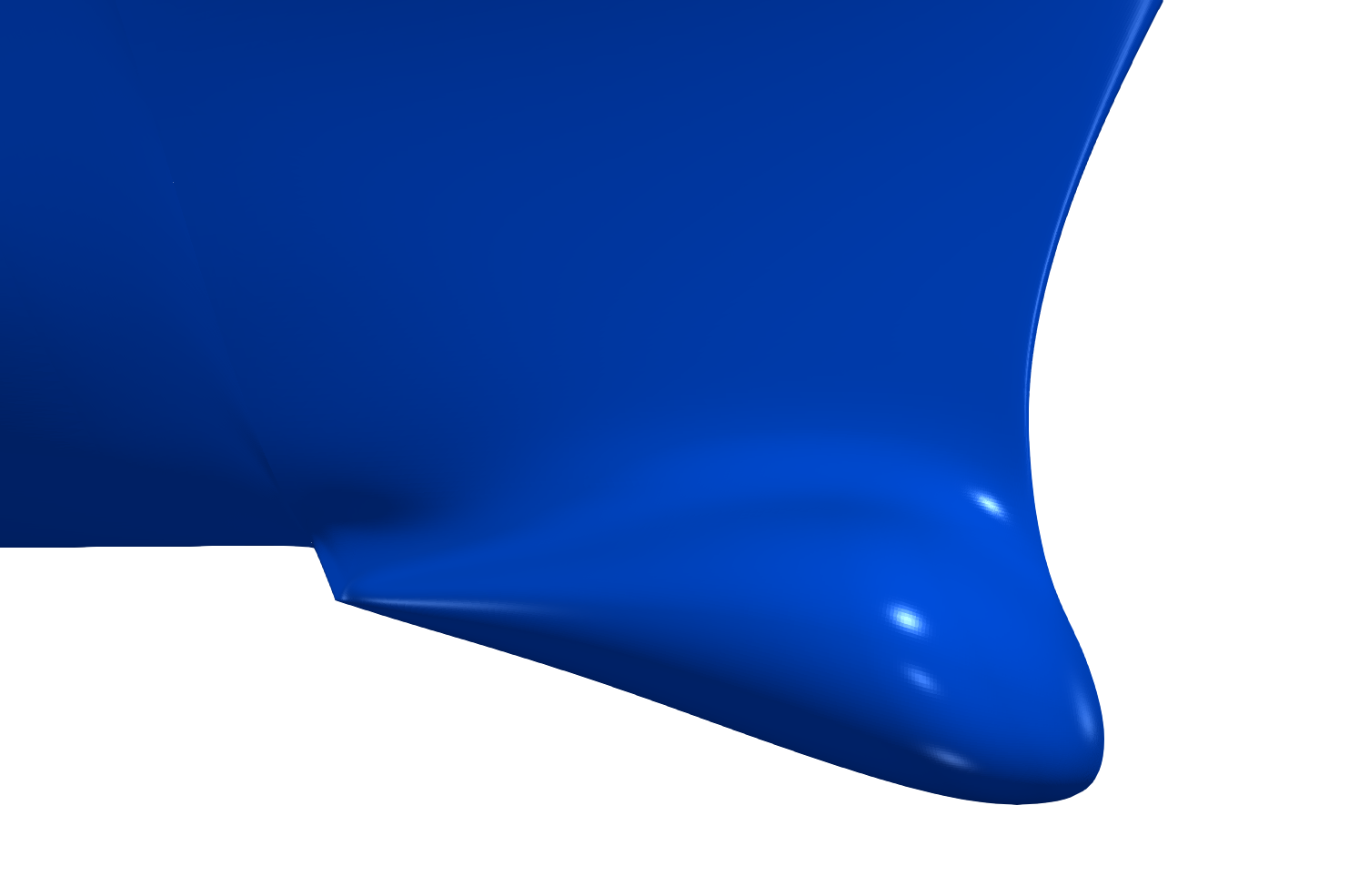}
\caption{Two different examples of deformations.}
\label{fig:deform_examples}
\end{figure}

In our case we have chosen $\boldsymbol{\mu} \in \mathbb{D} := [-0.3,
0.3]^8$, that results in a wide range of different bulbous
bow configurations (see Figure~\ref{fig:deform_examples} for an idea of
different possible deformations). For sake of clarity we underline
that the undeformed original domain is obtained setting all the
geometrical parameters to 0. The main purpose of this work is to
carry out a rather wide exploration of the shape parameter space, to
assess the dependence of the output parameters on the shape morphing parameters.
That is why the parameters bounds are chosen so as to be able to obtain somewhat
significant --- and yet physically meaningful --- deformations of the bow bulb.

\section{System evolution reconstruction with dynamic mode decomposition}

As mentioned, the unsteady and nonlinear fluid dynamic model adopted results in
rather expensive simulations. Of course this could in principle reduced through
full parallelization of the software developed which, however, has yet to be
carried out. In this work, the computational cost of each simulation carried
out has been reduced through the first application application of the dynamic
mode decomposition (DMD) to our fully nonlinear potential solver output.

Dynamic mode decomposition (DMD) is a data-driven algorithm that provides a
finite approximation of the infinite dimensional Koopman operator
(see~\cite{koopman1931hamiltonian}). Proposed
in~\cite{schmid2010dynamic} for fluid dynamics analysis, this
technique has become popular in the last years mainly because ({\it i}) it allows to
approximate nonlinear dynamics through low-rank structures that evolve in time
and ({\it ii}) it relies only on the data, avoiding assumptions on the
underlying system. We have implemented this algorithm, as well as many of its
variants, in an open source Python package on GitHub, called PyDMD~\cite{demo2018pydmd}. In this
section we will introduce the DMD algorithm and we will discuss the fluid dynamic
problem at hand as an example of its application. For an application
of DMD to snapshots obtained from the solution of RANS equations see~\cite{demo2018shape}.

\begin{figure}[p]
\centering
\includegraphics[width=.7\textwidth]{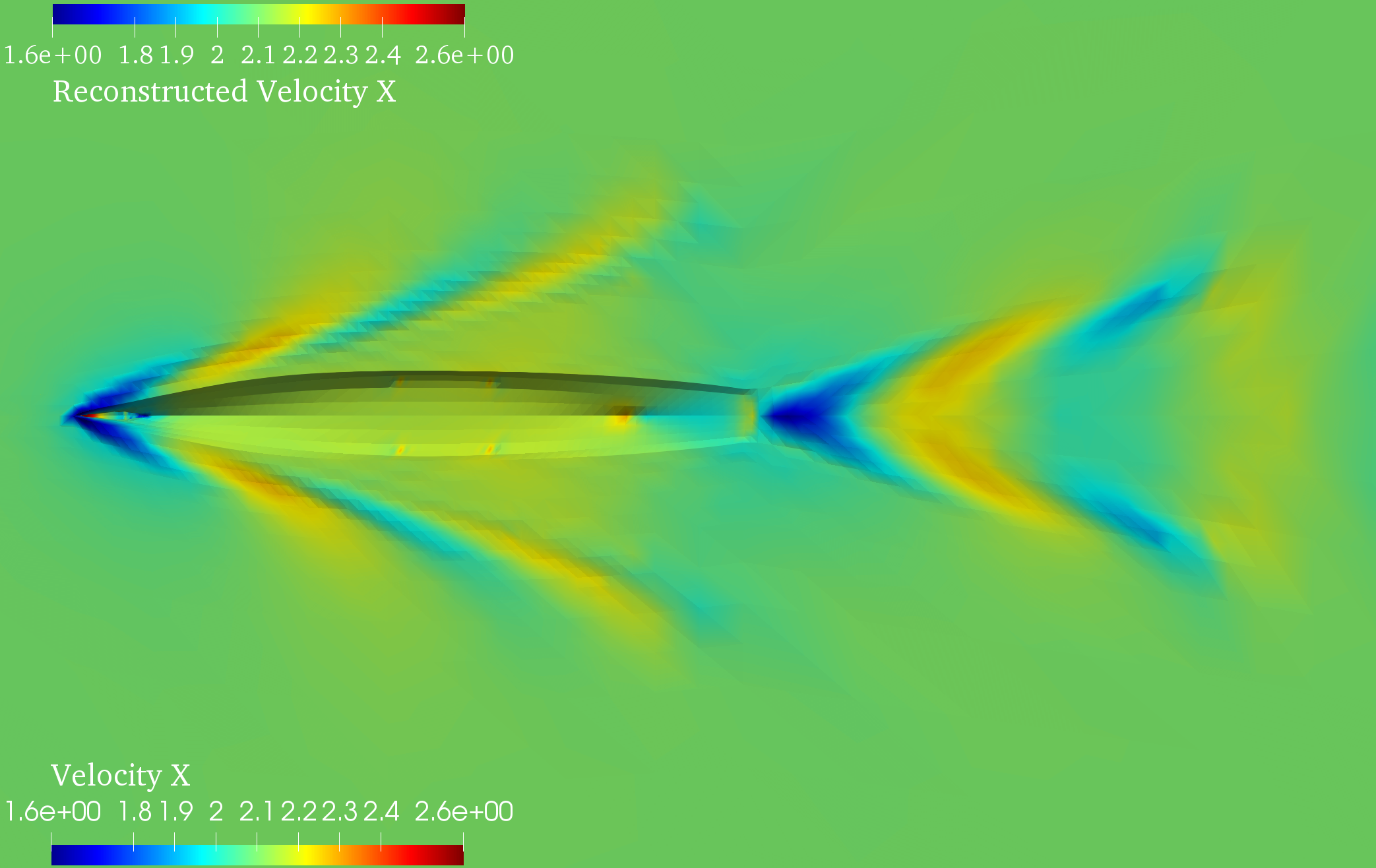}
\includegraphics[width=.7\textwidth]{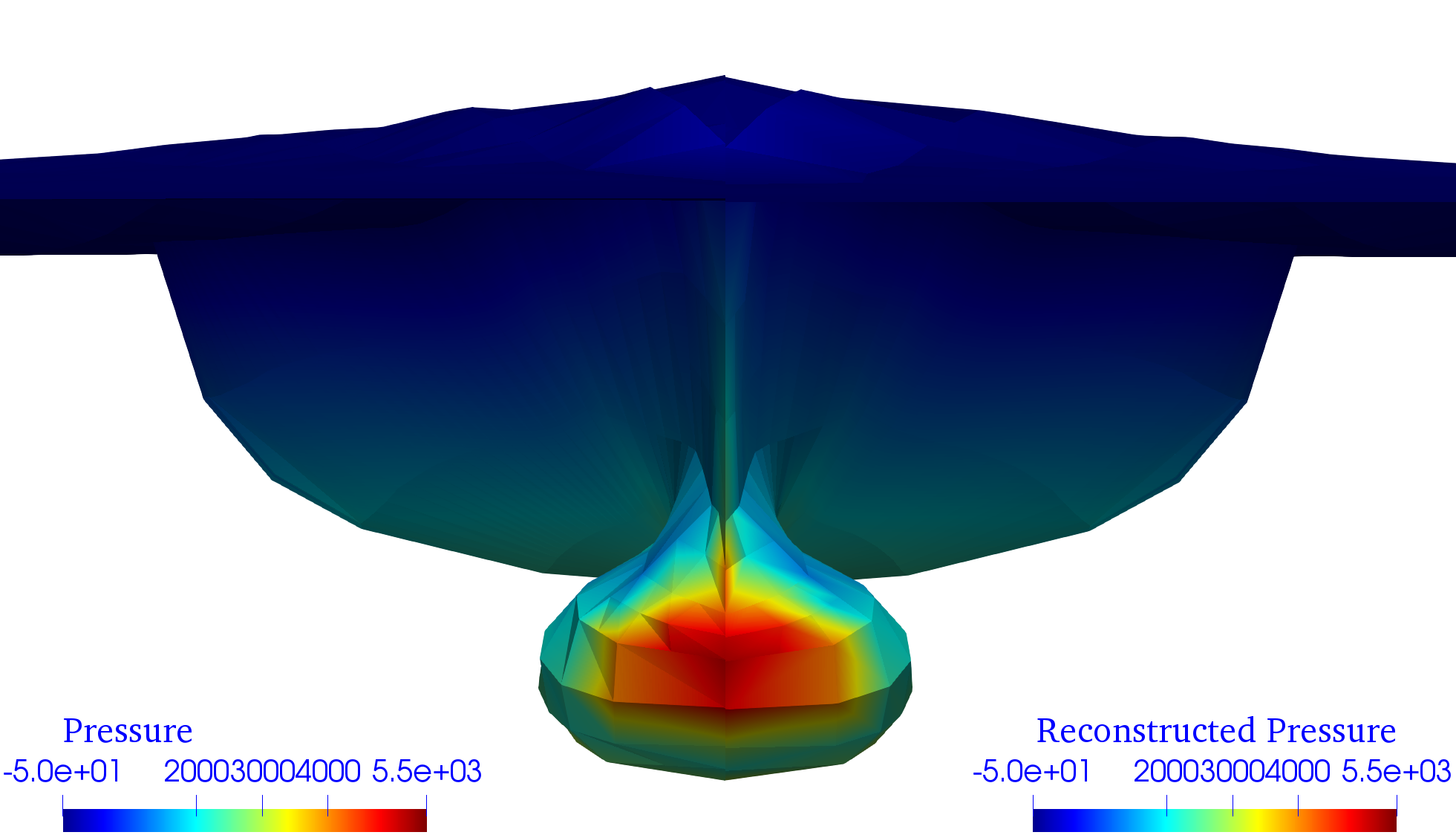}
\caption{A comparison between the original output of the fluid dynamic simulations and
         the DMD reconstructed one. The top image represents a contour plot of water $X$ velocity
         around a morphing of the DTMB 5415 hull advancing at $\text{Fr}=0.28$ in calm water.
         The top half of the plot depicts the reconstructed fluid velocity, while the bottom half
         represents the one resulting from the high-fidelity computation. The bottom image
         is a front view of a contour plot of the pressure field on the hull and water surface.
         In this case, the left part of the plot refers to the high-fidelity solution, while the
         right half refers to the reconstructed pressure.}
    \label{fig:dmdoutput}
\end{figure}
Let the variable $x_k$ represent the state of the evolving system at time $t_k = k\Delta
t$. Basically, we want to find a linear finite dimensional Koopman operator
$\mathbf{A}$ such that:
\begin{equation}
\label{eq:aevolution}
    x_{k+1} = \mathbf{A}x_k
\end{equation}
In order to build this operator, we collect a series of data vectors
$\{x_i\}_{i=1}^l$, which we will refer to as snapshots from now on, and which
represent the time-equispaced system states. We assume all the snapshots have the
same dimension, that is $x_k
\in \mathbb{R}^n$ for all $k = 1, \dotsc, l$, and we assume the dimension $n$ of a
snapshot is larger that the number of snapshots $l$, i.e.\ $n > l$.
We arrange the snapshots in two matrices, $\mathbf{S}$ and $\mathbf{\dot{S}}$, as
\begin{equation}
\label{eq:matarranged}
\mathbf{S} =
 \begin{bmatrix}
  x_1^1   & x_2^1  & \cdots & x_{l-1}^1 \\
  x_1^2   & x_2^2  & \cdots & x_{l-1}^2 \\
  \vdots  & \vdots & \ddots & \vdots    \\
  x_1^n   & x_2^n  & \cdots & x_{l-1}^n 
 \end{bmatrix},\quad\quad
    \mathbf{\dot{S}} =
 \begin{bmatrix}
  x_2^1   & x_3^1  & \cdots & x_l^1  \\

  x_2^2   & x_3^2  & \cdots & x_l^2  \\
  \vdots  & \vdots & \ddots & \vdots \\
  x_2^n   & x_3^n  & \cdots & x_l^n 
 \end{bmatrix}
\end{equation}
in order to build the linear operator by minimizing $\|\mathbf{\dot{S}} -
\mathbf{A}\mathbf{S}\|_2$. We underline that each column of $\mathbf{\dot{S}}$
contains the state vector at the next timestep of the one in the corresponding
$\mathbf{S}$ column.
Hence, the best-fit matrix $\mathbf{A}$ is given by:
\begin{equation}
\label{eq:adefinition}
    \mathbf{A} = \mathbf{\dot{S}}\mathbf{S}^\dagger
\end{equation}
where the symbol $^\dagger$ denotes the Moore-Penrose pseudo-inverse. Since the
snapshots usually have high dimension for complex systems, the matrix
$\mathbf{A}$ becomes very large and it is difficult to manipulate. The DMD
algorithm projects the data onto a low-rank subspace defined by the Proper Orthogonal
Decomposition (POD) modes,
then computes the low-dimensional operator $\tilde{\mathbf{A}}$. This operator
is used to reconstruct the leading nonzero eigenvalues and eigenvectors of the
full-dimensional operator~$\mathbf{A}$ without ever explicitly
computing~$\mathbf{A}$. Using the truncated singular value decomposition of
matrix $\mathbf{S} \approx \mathbf{U}_r \mathbf{\Sigma}_r \mathbf{V}^*_r
$ with rank $r$, we can build the low-rank linear operator
as:
\begin{equation}
\mathbf{\tilde{A}}  = \mathbf{U}_r^* \mathbf{A} \mathbf{U}_r
    = \mathbf{U}_r^* \mathbf{\dot{S}} \mathbf{S}^\dagger \mathbf{U}_r
     = \mathbf{U}_r^* \mathbf{\dot{S}} \mathbf{V}_r \mathbf{\Sigma}_r^{-1}
 	\mathbf{U}_r^* \mathbf{U}_r
     = \mathbf{U}_r^* \mathbf{\dot{S}} \mathbf{V}_r \mathbf{\Sigma}_r^{-1}.
\end{equation}
We can now reconstruct the eigenvectors and eigenvalues of the matrix
$\mathbf{A}$ using the eigendecomposition $\mathbf{\tilde{A}} \mathbf{W} =
\mathbf{W} \mathbf{\Lambda}$. In detail (see~\cite{tu2014dynamic}), the DMD
modes~$\mathbf{\Theta}$ can be computed by projecting the low-rank eigenvectors
on the high-dimensional space $\mathbf{\Phi} = \mathbf{U}_r \mathbf{W}$
(\textit{projected modes}) or computing the eigenvectors of $\mathbf{A}$ as
$\mathbf{\Theta} = \mathbf{\dot{S}}\mathbf{V}_r \mathbf{\Sigma}_r^{-1} \mathbf{W}$
(\textit{exact modes}). Moreover, the eigenvalues of $\mathbf{\tilde{A}}$
correspond to the nonzero eigenvalues of $\mathbf{A}$, and they contain the
growth/decay rate and the frequencies of the corresponding modes.  We recall
Eq.~\eqref{eq:aevolution} and underline that $\mathbf{A} = \mathbf{\Theta}
\mathbf{\Lambda} \mathbf{\Theta^\dagger}$. The generic snapshot $x_{k+1}$ can be
reconstructed by premultiplying the first snapshot $k$ times by the linear
operator, such that $x_{k+1} = \mathbf{A}^kx_1 = (\mathbf{\Theta}
\mathbf{\Lambda} \mathbf{\Theta^\dagger}\dotsc \mathbf{\Theta} \mathbf{\Lambda}
\mathbf{\Theta^\dagger})x_1$. Hence the state of the system can be approximated,
for any time $t_{k+1}$ as:
\begin{equation}
    x_{k+1} = \mathbf{\Theta} \mathbf{\Lambda}^k \mathbf{\Theta}^\dagger x_1
\end{equation}
where the vector $\mathbf{\Theta}^\dagger x_1$ is usually called
\textit{amplitudes}.

\begin{figure}[H]
\centering
\includegraphics[width=.7\textwidth]{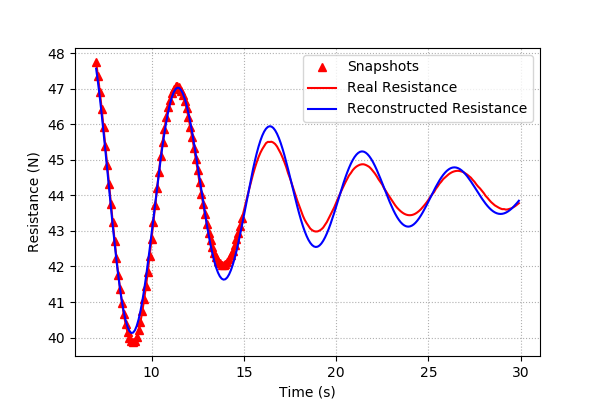}
\caption{Comparison of high-fidelity (red continuous curve) and DMD reconstructed (blue
         continuous curve) total resistance time history. The triangular markers in the plot
         denote the times at which the snapshots for the DMD have been stored.}
    \label{fig:dmdtrend}
\end{figure}
The application of DMD to the computational fluid dynamics simulations in this work
is carried out collecting the snapshots within the temporal window $t = [7\ s, 15\ s]$,
with $\Delta t = 0.1\ s$. As both the motion of the hull and of the free surface nodes
are computed in the simulations, the DMD algorithm is educated using snapshot vectors
which include the grid nodes coordinates, along with the flow velocity and pressure.
The DMD algorithm has been used to compute a rather accurate approximation of the whole
fluid dynamic problem solution up to $t = 30\ s$, at the mere cost of a post processing phase
requiring only few seconds. This led to a significant reduction of the computational cost,
dropping the time required for each $30\ s$ simulation from $24\ h$ to $10\ h$.

Figure~\ref{fig:dmdoutput} compares the high-fidelity fluid dynamic solution at $t = 30\ s$
with the one which has been reconstructed through the described DMD strategy.
The figure includes both pressure and velocity field plots, which have been split in half to
give a parallel view of both the high-fidelity and reconstructed field. Despite the
DMD has been extrapolated from snapshots that are alder than $t = 15\ s$, no difference
between original and reconstructed flow and water elevation fields is appreciable in the
contour plots. Figure~\ref{fig:dmdtrend} represents a comparison of time history
of the hull resistance 
computed from the high-fidelity and DMD reconstructed hull pressure field. Again, the
DMD approximated time history seems able to reconstruct with good approximation the
the high-fidelity values. This is especially true for the last instants of the simulation,
which are used to compute the final value of the total hull resistance, one of the output
parameters considered in this work.

\section{Parameter space verification and reduction by active subspaces}

The active subspaces (AS) property has been brought to attention recently through
the work of P.~Constantine~\cite{constantine2015active}. The AS property is a characteristic
of the scalar function relating the scalar output to the parameters $\mupar$, and
of a probability density function associated to such function.
By a qualitative standpoint, AS is typically exploited to assess whether the
parameter space allows for a significant --- and of course useful --- dimension reduction.
By a quantitative standpoint, it can also be used to assess the sensitivity of the
output with respect to each parameter considered. 
The main idea of AS is to operate in the parameters space, rescaling
the inputs $\mupar$ and then rotating them with respect to the
origin. In some cases (see~\cite{tezzele2017combined,tezzele2017dimension,tezzele2018model}), such procedure reveals in lower dimension
behavior of the output function $f(\mupar)$
(the total resistance or the hull sink and trim, in our case).
We underline that AS does not identify a subset of the
inputs as important, instead it identifies a set of important
directions in the space of all inputs. These directions (which are
linear combinations of the input variables) are the ones
along which the output function varies the most on average. When an active
subspace is identified for the problem of interest, it is
possible to perform different parameter studies.

Now we review how it is possible to find active subspaces. Let us
assume\footnote{In this section we will omit the dependence on
  $\mupar$. It should be understood that $f = f(\mupar)$, $\rho =
  \rho(\mupar)$, etc.} $f: \mathbb{R}^m \rightarrow \mathbb{R}$ is a
scalar function and $\rho: \mathbb{R}^m \rightarrow \mathbb{R}^+$ a
probability density function, where $m$ is the dimension of the
parameters. Since all the geometrical configurations can
be drawn with equal probability, a uniform probability
density will suffice in our case. In particular, we assume $f$
continuous and differentiable in the support of $\rho$, with
continuous and square-integrable (with respect to the measure induced
by $\rho$) derivatives. The active subspaces of the pair $(f, \rho)$
are the eigenspaces of the covariance matrix associated to the
gradients $\nabla_{\mupar} f$. This matrix, denoted by
$\boldsymbol{\Sigma}$, is the so-called
uncentered covariance matrix of the gradients of $f$ (among others
see~\cite{devore2015probability} for a more deep understanding of
these operators). Its elements are the average products of partial
derivatives of $f$, that is:
\begin{equation}
\label{eq:covariance}
\boldsymbol{\Sigma} = \mathbb{E}\, [\nabla_{\mupar} f \, \nabla_{\mupar} f
^T] = \int_{\mathbb{D}} (\nabla_{\mupar} f) ( \nabla_{\mupar} f )^T
\rho \, d \mupar ,
\end{equation}
where $\mathbb{E}[\cdot]$ is the expected value. We use a Monte Carlo
method to approximate the eigenpairs of $\boldsymbol{\Sigma}$ as in~\cite{constantine2015computing}:
\begin{equation}
\label{eq:covariance_approx}
\boldsymbol{\Sigma} \approx \frac{1}{N_{\text{train}}^{\text{AS}}} \sum_{i=1}^{N_{\text{train}}^{\text{AS}}} \nabla_{\mupar} f_i \,
\nabla_{\mupar} f^T_i ,
\end{equation}
where we draw $N_{\text{train}}^{\text{AS}}$ independent samples $\mupar^{(i)}$ from the
measure $\rho$ and where $\nabla_{\mupar} f_i = \nabla_{\mupar}
f(\mupar^{(i)})$.
The matrix $\boldsymbol{\Sigma}$ admits a real eigenvalue
decomposition because it is symmetric positive semidefinite, so we
have
\begin{equation}
\label{eq:decomposition}
\boldsymbol{\Sigma} = \mathbf{W} \mathbf{\Lambda} \mathbf{W}^T ,
\end{equation}
where $\mathbf{W}$ contains the eigenvectors and is in $O(m)$, the orthogonal group,
while $\mathbf{\Lambda}$ is the diagonal matrix of non-negative
eigenvalues arranged in descending order.

The lower dimensional parameter subspace is formed by selecting the
first $M < m$ eigenvectors. We underline that perturbations
in the first set of coordinates change $f$, on average, more than perturbations in
the second set of coordinates. We can discard the vectors
corresponding to the low eigenvalues since they are in the nullspace
of the covariance matrix. Doing so, we are able to construct an approximation
of $f$. To be more clear, let us partition $\mathbf{\Lambda}$ and
$\mathbf{W}$ as follows: 
\[
\mathbf{\Lambda} =   \begin{bmatrix} \mathbf{\Lambda}_1 & \\
                                     &
                                     \mathbf{\Lambda}_2\end{bmatrix},
\qquad
\mathbf{W} = \left [ \mathbf{W}_1 \quad \mathbf{W}_2 \right ],
\]
where $\mathbf{\Lambda}_1 = \text{diag}(\lambda_1, \dots, \lambda_M)$, and
$\mathbf{W}_1$ contains the first $M$ eigenvectors. The active subspace
is the the range of $\mathbf{W}_1$. 
We call inactive subspace the range of the remaining eigenvectors in
$\mathbf{W}_2$.
The linear combinations of the input
parameters with weights from the important eigenvectors are the active variables. 
By projecting the full parameter space onto the active
subspace we can approximate the behaviour of $f$. In particular we
have the following formulas for the active variable $\mupar_M$ and the
inactive variable $\etapar$:
\begin{equation}
\label{eq:active_var}
\mupar_M = \mathbf{W}_1^T\mupar \in \mathbb{R}^M, \qquad
\etapar = \mathbf{W}_2^T \mupar \in \mathbb{R}^{m - M} .
\end{equation}
Using Eq.~\eqref{eq:active_var} and the fact that $\mathbf{W} \in
O(m)$ we can express any point in the parameter space $\mupar \in
\mathbb{R}^m$ in terms of $\mupar_M$ and $\etapar$ as follows:
\[
\mupar = \mathbf{W}\mathbf{W}^T\mupar =
\mathbf{W}_1\mathbf{W}_1^T\mupar +
\mathbf{W}_2\mathbf{W}_2^T\mupar = \mathbf{W}_1 \mupar_M +
\mathbf{W}_2 \etapar.
\]
So it is possible to rewrite $f$ as
\[ 
f (\mupar) =  f (\mathbf{W}_1 \mupar_M + \mathbf{W}_2 \etapar) ,
\]
and, using only the active variables, we can construct a surrogate
quantity of interest $g$
\[
f (\mupar) \approx g (\mathbf{W}_1^T \mupar) = g(\mupar_M).
\]
In our pipeline, the surrogate quantity of interest $g$ will be
obtained by a response surface method.

\section{Numerical Results}
In this section we present the numerical results obtained by applying
all the methodologies presented above to the DTMB 5415 model hull.
The FFD morphing methodology illustrated has been used to generate
130 different deformations of the original hull. The parametrised shapes
correspond to uniform sampling points in the parameter space box
$\mathbb{D} = [-0.3, 0.3]^8$. Each IGES geometry produced has been then
used as the input of a high fidelity simulation in which the hull has been
set to advance in calm water at a constant speed corresponding to
$\text{Fr}=0.28$. Each high fidelity computation has been carried out to
simulate $15\ s$  of the flow past the hull after it has been
impulsively started from rest. Between the 7th and 15th second of the
high-fidelity simulations, the solver saved the full flow field
at sampling intervals $\Delta t = 0.1\ s$. Such flow field snapshots have
been used to feed the DMD algorithm implemented, and complete the fluid
dynamic simulations until convergence to the regime solution was reached
at $t=30\ s$. The reconstructed flow fields have been finally used to
evaluate the hull total resistance and the hydrodynamic trim position,
which are the output performance parameters considered in this work.
The dataset composed by the output of the simulations has been divided in a
train dataset (75\% of the outputs) used to train the AS algorithm and
a test dataset (25\% of the outputs) used to validate the methodology.

\begin{figure}[h]
\centering
\includegraphics[width=.46\textwidth]{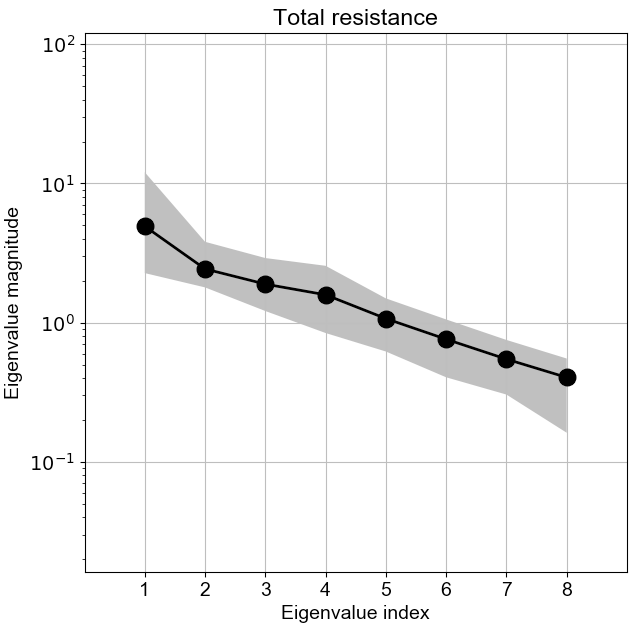}
\includegraphics[width=.48\textwidth]{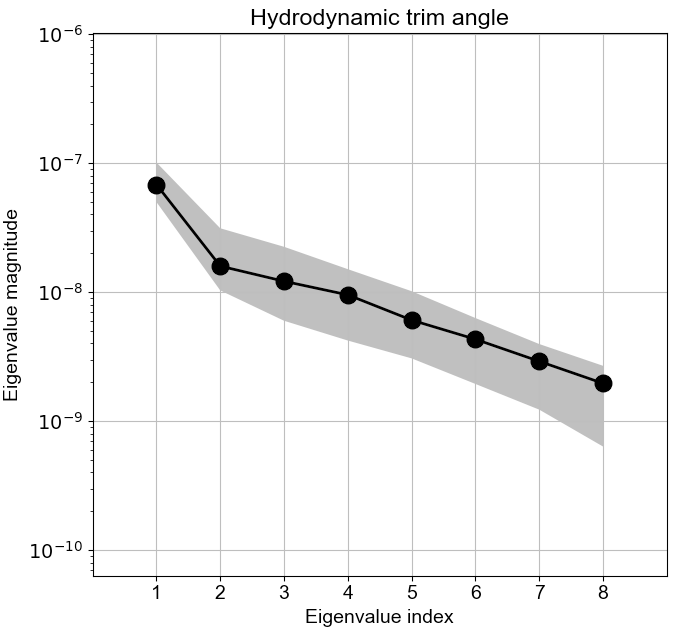}
\caption{Eigenvalues estimates of the matrix $\boldsymbol{\Sigma}$ for
         the total resistance (on the left) and the hydrodynamic trim angle (on the right).
         The black dots in the plot indicate the eigenvalues, which the grey area
         is defined by the bootstrap intervals.}
\label{fig:eigs}
\end{figure}

Figure~\ref{fig:eigs} depicts the eigenvalues estimates
of the matrix $\boldsymbol{\Sigma}$  (black dots and line) and also
displays the bootstrap intervals (corresponding to the grey area surrounding
the eigenvalues lines). In the figure, the left plot is referred to
the total resistance, while the right one is displays the hydrodynamic
trim angle.

\begin{figure}[h]
\centering
\includegraphics[trim={3cm 0 0 1cm},clip,width=.48\textwidth]{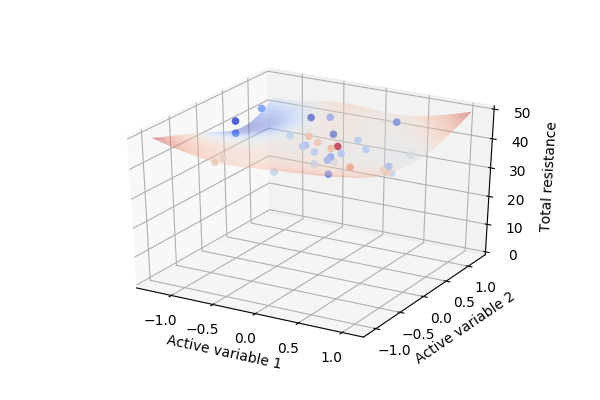}
\includegraphics[trim={3cm 0 0 1cm},clip,width=.48\textwidth]{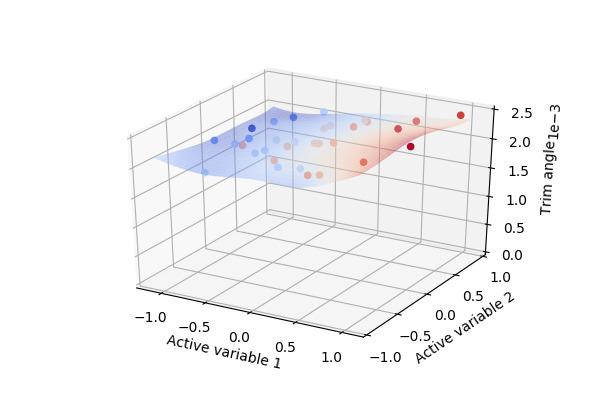}
\caption{Comparison between the quartic surface response obtained with two active variables
         and the true output for the test dataset, indicated by the dots. The left plot refers
         to the total resistance output, while the right one to the hydrodynamic trim angle.}
\label{fig:ss2}
\end{figure}

The plots indicate that a factor of at least 10 exists between the
highest and lowest $\boldsymbol{\Sigma}$ eigenvalues. Such difference
is clearly more pronounced when the hydrodynamic trim is the output considered.
Yet, the plots also show that the eigenvalues magnitude is rather evenly distributed
across the range they span. More precisely, the absence of a major gap between
the higher module eigenvalues and the lower module ones, is suggesting that
the active subspace is most revealing a clear cut low
dimensional behaviour of the target functions with respect to the
active variables, as is the case for different applications~\cite{tezzele2017combined,tezzele2017dimension}.
Yet, especially in the case of the hydrodynamic trim angle output, the
first eigenvalue module is considerably higher than that of the remaining ones. That
is why, for the hydrodynamic trim it was possible to compute a bivariate surface response using
the first two active variables, obtained as linear combinations of the original parameters
with coefficients obtained by the eigenvectors corresponding to the two highest module eigenvalues.
Figure~\ref{fig:ss2} shows the quartic surface that best approximates
the training dataset in the sense of least squares, along with the
points in the test dataset (which are indicated by the dots). Each
point represents the value of the target function $f(\mupar)$ against
the active variables $\mupar_M = \boldsymbol{W}_1^T \mupar \in \mathbb{R}^2$.
As can be appreciated, the points corresponding to the true output are not distributed
randomly in the space, but tend to be somewhat clustered around the surface. This is particularly true
when the output parameter considered is the hydrodynamic trim angle (right plot) for which, as we have seen,
the $\boldsymbol{\Sigma}$ eigenvalues corresponding to the first active variable was significantly
higher then the remaining ones. Thus, whenever the gap between the leading $\boldsymbol{\Sigma}$
eigenvalues allows for it, the AS algorithm is able to successfully identify a set active variables
upon which the output is --- with reasonable approximation --- exclusively depending.

To provide a more quantitative assessment of how much such approximation is in fact reasonable, and
depends on the gap existing between the first $\boldsymbol{\Sigma}$ eigenvalues and the
following ones, we computed the average error as the average, among 10 different
eigenvalues estimates, of the root mean square error divided by the maximum range of variation
of the target function. As expected, the accuracy of the two dimensional surface response
predictions for the hull total resistance is rather low, and only a 20\% average error is
obtained, that is about $1.8$ Newton. The average error obtained for the hydrodynamic trim angle output is approximatively
11\%, that corresponds to an absolute rmse of $8\cdot 10^{-5}$, confirming that such output can be better represented through AS. Thus, for all those
applications for which 10\% can be considered an acceptable error --- as might be the case
for early design stages --- AS provides a recipe to reduce the parameter space from 8 to 2
variables. In addition, the error analysis provides
a further confirmation of the fact that once the eigenvalue analysis is carried out, the
detection of a possible cliff in the eigenvalue curve (see Figure~\ref{fig:eigs}) is a measure
of how well AS will perform. So, since the computational cost of the post processing operations
required for the AS algorithm proposed is marginal with respect to that of the high-fidelity
simulations, it should be always worth checking if AS could be used to obtain a significant
drop in the parameter space dimension.

\section{Conclusions and further developments}

This work presented an application to naval hydrodynamics
simulations of the active subspaces algorithm. Such
algorithm can lead to a significant reduction of the parameter space
dimension. The test case considered consisted in the flow
past several parametrised modifications of the shape of a DTMB 5415 model hull. Such
modified shapes were obtained by means of the free form deformation involving 8 parameters,
and used for fully nonlinear potential flow simulations. The computational
cost of the whole simulation campaign has been significantly reduced through
the application of dynamic mode decomposition strategy to speedup
the convergence of the unsteady simulation to the finial regime solution. 
The contribution describes the whole pipeline composed of shape parametrization,
hydrodynamic simulations, and data post processing, all based on in house
developed software.

The results obtained confirm that AS can be conveniently used to reduce
the parameter space dimension when a significant gap is detected between the highest and
lowest eigenvalues of the uncentered covariance matrix of the output gradients
with respect to the input parameters. Yet, in the present case the hull total
resistance (one of the two outputs considered) did not present such feature,
resulting in a diminished effect of the AS parameter dimension reduction. As
for the hydrodynamic trim angle, AS is able to provide more accurate two dimensional
approximation of the output dependence on the eight shape parameters considered.
The reduced computational cost of the post processing steps required for the AS
analysis suggests that it can be used in any contest to check whether a
reduction of the parametric space dimension is convenient.

Future work will focus on several areas to improve physical and mathematical
aspects of the algorithms presented. Further developments of the whole pipeline
will be implemented to better integrate its different parts and automate the
simulation campaign and post processing processes.

\section*{Acknowledgements}
This work was partially performed in the context of the project SOPHYA,
``Seakeeping Of Planing Hull YAchts'', supported by Regione
FVG, POR-FESR 2014-2020, Piano Operativo Regionale Fondo Europeo per lo Sviluppo Regionale,
partially funded by the project HEaD, ``Higher
Education and Development'', supported by Regione FVG - European
Social Fund FSE 2014-2020, and by
European Union Funding for Research and Innovation --- Horizon 2020 Program ---
in the framework of European Research Council Executive Agency: H2020 ERC CoG
2015 AROMA-CFD project 681447 ``Advanced Reduced Order Methods with
Applications in Computational Fluid Dynamics'' P.I. Gianluigi Rozza,
as well as MIUR FARE-X-AROMA-CFD project.


\end{document}